**ON PROBLEM OF CONFLICT INTERACTION OF TWO PLAYERS IN THE NETWORK ENVIRONMENT**

*Oleksii Ignatenko*

Institute of software systems NAS of Ukraine
40 Glushkov Prsp., 03187 Kyiv, Ukraine,
tel.: +380 526 60 25, E-Mail: oignat@isofts.kiev.ua

**Abstract.** In developing the network management systems often encountered a situation where there is a need to deal with malicious intrusion attempts aimed at obstructing the work. In this case, methods of management based on fluid models and a discrete controlled random walk model requires certain changes. In this work a generalization of the usual models of single server queue which can describe the attacking action and behaviour of the system of protection is presented. Formulation and analysis of the problem are proposed. There are conditions found under which the game can be finished.

## Introduction

Areas where one can apply network models are very different – information networks, telecommunications, gas transportation and energy systems, distributed production processes. Information network applications are especially important. Today networks are related to almost all sides of human activity. As a result, security and reliability of information flows directly affect the quality of service, efficiency and overall economic development of entire industries. Networks complexity and interconnections are growing constantly. Today most of organizations depend on reliable working of the information networks. On the other hand, the vulnerability of networks is also increasing every year. This is due mainly to the fact that more potential attackers have access to resources with networks growing.

Nowdays, network is a complex (often nonlinear) system with significant unknown disturbance. Our influence possibilities are constrained by some parameters called control parameters. An important point of control theory is its flexible approach to modeling. We typically consider a linear, deterministic model. In the work [4] said an idea that the control system should be simplest model that capture essential feature of the system to be controlled.

Let us define network terminology. The network here consists of a finite set of nodes. Each node contains a set of buffers. These buffers can hold finite amount of packets. One or more servers process packets at a given buffer.

A general stochastic model defined in the work [4]. Let us make a short summary.

Define the vector valued queue-length process $Q$ that evolves on $R_+^N$. The vector valued cumulative allocation process $Z \in R_+^{N_u}$. The *i*th component of $Z(t)$, denoted $Z_i(t)$, is equal to the cumulative time that the activity *i* has run up to time *t*. The evolution of the queue-length process is described by the vector equation,

$$Q(t) = Q(t_0) + B(Z(t)) + A(t)\,,\ t \geq 0\,. \tag{1}$$

where the process $A$ may denote a combination of exogenous arrivals to the network, and exogenous demands for materials from the network. The function $B$ represents the effects of (possibly random) routing and service rates.

The cumulative allocation process and queue-length process are subject to several hard constraints:

(i) The queue-length process is subject to the state space constraint,

$$Q(t) \in X\,,\ t \geq 0\,,\ X \subset R_+^N\,. \tag{2}$$

(ii) The control rates are subject to linear constraints,

$$C(Z(t_1) - Z(t_0)) \leq (t_1 - t_0) \cdot e\,,\ Z(t_1) - Z(t_0) \geq 0\,,\ 0 \leq t_0 \leq t_1\,, \tag{3}$$

where the constituency matrix $C$ is a $N_m \times N_u$ matrix, $e = (1,...,1) \in R^{N_m}$.

Stochastic models such as (1) have been by far the most popular in queuing theory. An idealization is the linear fluid model, described by the purely deterministic equation,

$$q(t) = q(t_0) + Bz(t) + \alpha(t - t_0)\,,\ t \geq 0\,,\ q(t) \in R_+^N\,. \tag{4}$$

where the vector $q(t)$ evolves in the state space $R_+^N$ and the (cumulative) allocation process $z(t)$ evolves in $R_+^{N_u}$. We again assume that $z(0) = 0$, and for each $0 \leq t_0 \leq t_1$

$$C[z(t_1) - z(t_0)] \leq (t_1 - t_0) \cdot e\,. \tag{5}$$

$$z(t_1) - z(t_0) \geq 0\,.$$

The fluid model can also be described by the differential equation

$$\frac{dq}{dt} = Bu + \alpha\,, \tag{6}$$

where control $u(t)$ a is measurable function such that $Cu \leq e$, $u(t) \geq 0$. We choose a control function $u(t)$ in purpose of minimizing vector $q(t)$ under several special conditions (e.g. for minimal time).

Solution of (6) can be found using following formula:

$$q(t) = q(t_0) + \int_{t_0}^{t} Bu(\tau)d\tau + \alpha(t - t_0)\,.$$

Taking into account $z(t) = \int_{t_0}^{t} u(\tau)d\tau$ we obtain that expression above is equal to (4). Two different symbols are used to denote the state processes for the stochastic and fluid models since much of the development to follow is based on the relationship between the two models. In particular, the fluid model can be interpreted as the mean flow of the stochastic model (1) on writing,

$$Q(t) = Q(t_0) + A(t) + B(Z(t)) = Q(t_0) + BZ(t) + \alpha(t - t_0) + N(t)\,,\ t \geq 0\,, \qquad (7)$$

where $\alpha$ and $B$ are interpreted as average values of $A$ and $B$, and

$$N(t) = [A(t) - \alpha(t - t_0)] + [B(Z(t)) - BZ(t)]\,.$$

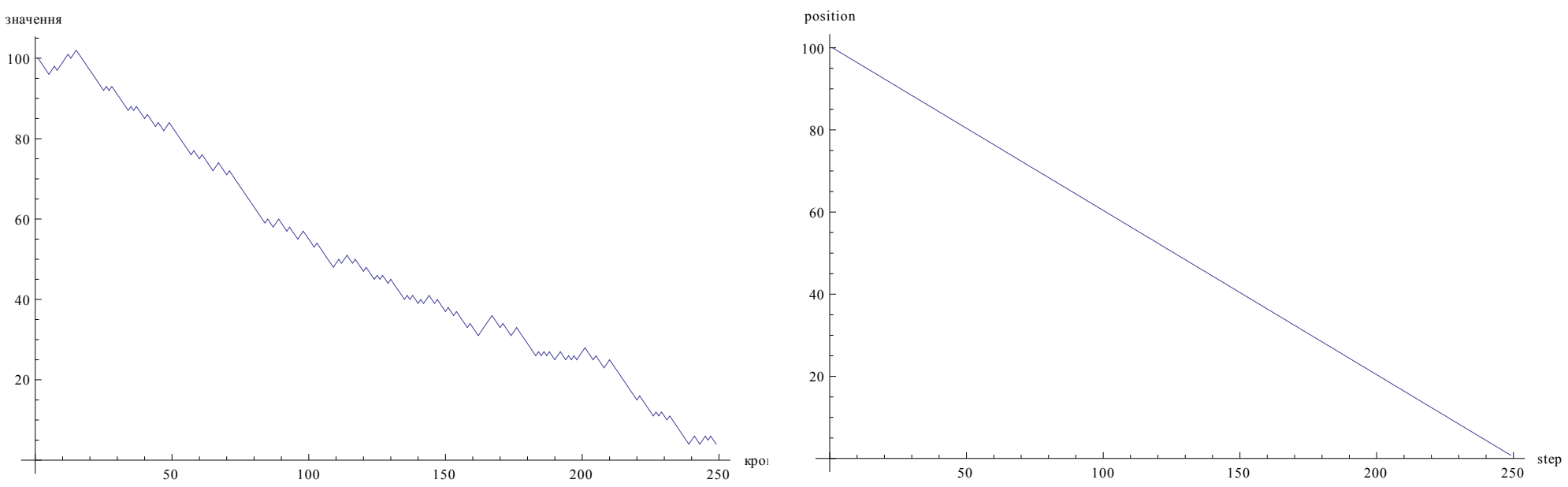

Figure 1. Stochastic and fluid model.

Typical assumptions on the stochastic model (7) imply that the mean of the process $\{N(t)\}$ is bounded as a function of time, and its variance grows linearly with $t$. Under these conditions (1) can be loosely interpreted as a fluid model subject to the additive disturbance $N(t)$. Example of behavior is presented on Fig. 1.

The theory of conflict controlled processes covers one of the divisions of mathematical theory of control which examines objects dynamic under conflict and uncertainty conditions. This conditions or factors may be different (technical problems, network topology bottlenecks, malicious bots) but finally they try to interrupt normal system's work.

Consider a conflict controlled object whose dynamics is described by the linear differential equation:

$$\frac{dq(t)}{dt} = Aq(t) + Bu(t) - Cv(t) + \alpha(t) \qquad (8)$$

where $q(t) \in R_+^n$ is an vector of an object state, $u(t) \in U \subset R^m$ is a parameter of first player, $v(t) \in V \subset R^l$ is a parameter of second player, $\alpha(\cdot): R \to R^n$ is a vector function of packets arriving. $A, B, C$ – real matrixes of orders $n \times n$, $n \times m$ i $n \times l$, respectively. A game (8) is played as follows: at moment of time $t \geq 0$ players set their control functions $v(t)$ and $u(t)$. Game is finished when $q(t) = 0$. First player tries to finish game as soon as possible. On the contrary, second player tries to counteract with all his recourses. If $\alpha(t) \equiv 0$ then equation (8) turns to classic conflict-controlled process [6]

$$\frac{dq(t)}{dt} = Aq(t) + Bu(t) - Cv(t) \qquad (9)$$

Problems of conflict processes of two players (9) are investigated in many works (e.g. [6 – 8]). In this work we focus on the first direct method of Pontryagin [7, 8].

## 2. Linear game model of the single server queue

The single server queue is a useful model for a dynamic investigation of very different systems. In this chapter we introduce an extension of the fluid model [4] in the case of conflict process. Let us consider dynamic system (Fig. 2). This system defined for an initial condition $q(t) = (q_1(t), q_2(t)) \in R^2$ by the system of linear equations:

$$\dot{q}_1(t) = \alpha(t) + k \cdot q_2(t) - u_1(t),$$
$$\dot{q}_2(t) = v(t) - u_2(t),\ t \geq 0. \qquad (10)$$

Phase state is described by the phase vector $q(t) = (q_1(t), q_2(t)) \in R^2$, where $q_1(t) \in R_+$ is the queue length at time $t$.

The queue length is subject to the linear phase constraint

$$0 \leq q_1(t) \leq q_1^{\max} \text{ for all } t \geq 0.$$

Parameter $q_2(t) \in R_+$ is associated with the attacker player which trying to overwhelm our network (or maximize $q_1(t)$ in other words) using control parameter $v(t)$, $v(t) \geq 0$, $v(t) \leq \nu$. By choosing $v(t)$ at time $t$ attacker sets attack power ($q_2(t)$) which is subject to the linear phase constraint

$$0 \leq q_2(t) \text{ for all } t \geq 0.$$

Parameter $q_2(t)$ influences on queue $q_1(t)$ with a coefficient $k \geq 0$.

The other player – defender – has a control vector $u(t) = (u_1(t), u_2(t))$. He divides his control resources between two directions: $u_1(t)$ - for service of the arrived packets, $u_2(t)$ - for counteraction of the attacker activity.

Control parameter $u(t)$ is subject of following constraints

$$u_1(t) \geq 0,\ u_2(t) \geq 0,\ u_1(t) + u_2(t) \leq \mu, \text{ for all } t \geq 0.$$

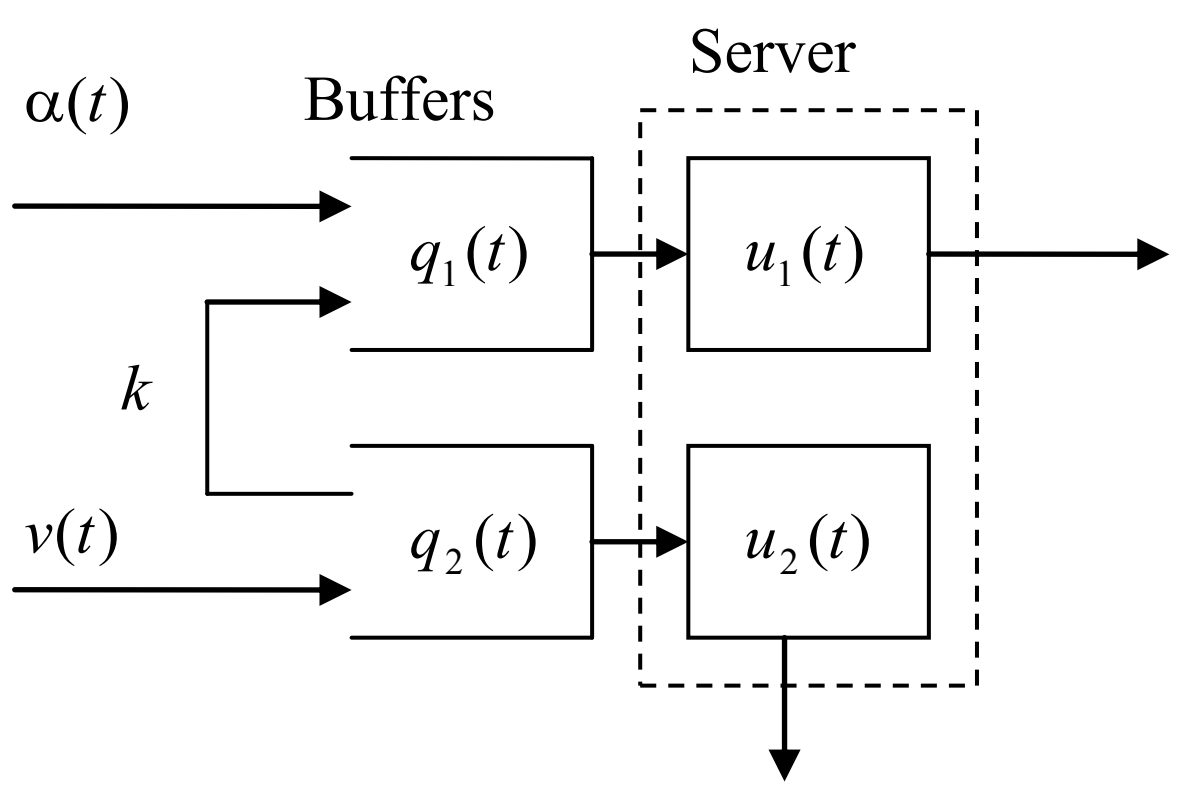

Figure 2. Single server queue with counteraction

Suppose that defender has information about $\alpha(t)$, $v(t)$ and $q(t)$ at time $t$.

Function $\alpha(\cdot): R \to R^n$ describe packets arrival service at time $t$. Suppose the following assumptions hold:

(i) $0 \le \alpha(t) \le \alpha_{\max}$ for all $t \ge 0$;

(ii) $\alpha(t)$ - continuous function.

If $q_2(0) = 0$ and $v(t) \equiv 0$, $\alpha(t) \equiv \alpha_{\max}$, then equations (10) are famous single server fluid model [].

Let us recall several results about the fluid model. The fluid model is called stabilizable if, for any $q_1(0)$, there exists strategy $u_1(t)$ and time $T \ge 0$ such that

$$q_1(T) = q_1(0) + \int_0^T (\alpha_{\max} - u_1(\tau))d\tau = 0. \quad (11)$$

The minimal draining time, denoted $T^*(q_1(0))$, is defined for $q_1(0) \in R$ as the smallest $T \ge 0$ satisfying (11).

For two states $q_1(t_1)$ and $q_1(t_2)$ the minimal time to reach $q_1(t_2)$ from $q_1(t_1)$ over all strategies is denoted $T(q_1(t_1), q_1(t_2))$.

This model is stabilizable if and only if $\frac{\alpha_{\max}}{\mu} < 1$. Time optimal strategy gives by following formula

$$u(t) = (u_1(t), 0), \text{ where } u_1(t) = \begin{cases} \mu, & q_1(t) > 0 \\ 0, & q_1(t) = 0 \end{cases}. \quad (12)$$

This is so-called non-idling policy. Minimal draining time in this case is

$$T^*(q_1(0)) = \frac{q_1(0)}{\mu - \alpha_{\max}}.$$

Our task is to obtain similar result for the equation (10). More formally, we should find conditions for $q_0$, such that following statement holds:

(i) $T^*(q_0) < +\infty$.

(ii) There are exists an admissible strategy $u(t, q(t), v(t))$ such that $q(T^*(q_0)) = 0$, where $q(t)$ is the solution of (10). We will call this strategy $u(t, q(t), v(t))$ - the solution of game (10).

Note that this result must be achieved over all an admissible functions $v(t)$.

Let us write equations (10) in more general form. Denote $A = \begin{pmatrix} 0 & k \\ 0 & 0 \end{pmatrix}$, then

$$\dot{q}(t) = Aq(t) + \begin{pmatrix} \alpha(t) \\ v(t) \end{pmatrix} - \begin{pmatrix} u_1(t) \\ u_2(t) \end{pmatrix} \quad (13).$$

Control sets are respectively $U = \{u \in R_+^2 : u_1 + u_2 \le \mu\}$ and $V = \{v \in R_+^2 : v_1 \le \alpha_{\max}, v_2 \le \nu\}$. Terminal set is $M = \{(0,0)\}$. Let us solve this problem using idea of the first direct method of Pontryagin [7, 8]. Let us consider the *Pontryagin map* $\omega(t) = \bigcap_{v \in V} \left( e^{At}U - e^{At}v \right)$.

$$\bigcap_{v \in V} \left( U - v \right) = \{x \in R^2 : x_1 + x_2 \le \mu - v - \alpha_{\max}\} =$$

$$= co\{(0,0), (\mu - v - \alpha_{\max}, 0), (0, \mu - v - \alpha_{\max})\} .$$

Since $e^{At}$ is linear operator, we see that

$$\bigcap_{v \in V} \left( e^{At}U - e^{At}v \right) = e^{At} \bigcap_{v \in V} \left( U - v \right).$$

$$\omega(t) = co\{(0,0), (\mu - \nu - \alpha_{\max}, 0), (kt(\mu - \nu - \alpha_{\max}), \mu - \nu - \alpha_{\max})\} .$$

**Condition 1.** Mapping $\omega(t) \neq \varnothing$ for all $t \ge 0$.

Condition 1 holds for all $t \ge 0$, if $\mu > \nu + \alpha_{\max}$.

As shown in [7] by M. Nikolsky (for game (10) without phase constraints) if condition 1 holds and exists time $T \ge 0$ such that $e^{AT} q_0 \in \int_0^T \omega(\tau) d\tau$ then the solution $u(t)$ could be constructed.

**Condition 2.** For initial state $(q_1(0), q_2(0))$ a following inequality holds

$$q_1^{\max} - q_1(0) \ge \frac{k}{2(\mu - \nu - \alpha_{\max})} \left( q_2(0) \right)^2 .$$

**Theorem.** Consider the game (10) satisfying conditions 1,2. Then, there are exists the solution $u(t)$ such that $q(T^*(q_0)) = 0$.

**Proof.**

Without loss of generality $q(0) \neq (0,0)$. Let us define moments of time

$$T_1 = \min\{ t > 0 : q_1(t) = 0\} ;$$

$$T_2 = \min\{ t \ge 0 : q_2(t) = 0\} .$$

Then, define $u(t)$

$$u(t) = \begin{cases} (\alpha(t), \mu - \alpha(t)), & t \in [0, T_2] \\ (\mu - v(t), v(t)), & t \in [T_2, T_1] \end{cases} \quad (14)$$

First we show then the strategy $u(t)$ is admissible. The strategy $u(t)$ is defined for all $t \geq 0$ and $u_1(t) + u_2(t) = \mu$.

Since condition 1 holds, it follows that $u_1(t) \geq 0$, $u_2(t) \geq 0$ for all $t \geq 0$.

Substituting (14) into (10) we obtain

For $t \in [0, T_2]$

$$\dot{q}_1(t) = k \cdot q_2(t),$$

$$\dot{q}_2(t) = v(t) + \alpha(t) - \mu.$$

For $t \in [T_2, T_1]$

$$\dot{q}_1(t) = v(t) + \alpha(t) - \mu,$$

$$\dot{q}_2(t) = 0.$$

Denote $\mu - \alpha_{\max} - \nu$ as $\varepsilon$, then

$$\dot{q}_2(t) \leq -\varepsilon,$$

$$q_2(t) \leq q_2(0) - \varepsilon t.$$

From the last statement follows that $q_2(t) = 0$ for moment of time $t \leq T_2^* = \dfrac{q_2(0)}{\varepsilon}$. Subsequently, we obtain that $T_2 \leq T_2^*$.

Similarly:

$$q_1(t) = q_1(0) + \int_0^t k q_2(\tau) d\tau$$

$$q_1(t) = q_1(0) + k q_2(0) t - \int_0^t k(v(\tau) + \alpha(\tau) - \mu) d\tau$$

$$q_1(t) \leq q_1(0) + k q_2(0) t - k\varepsilon \frac{t^2}{2}$$

$$q_1(T_2) \leq q_1(0) + k \frac{(q_2(0))^2}{\varepsilon} - k\varepsilon \frac{1}{2} \left( \frac{q_2(0)}{\varepsilon} \right)^2$$

$$q_1(T_2) \leq q_1(0) + \frac{k}{2} \frac{(q_2(0))^2}{\varepsilon}$$

$$q_1(T_2) \leq q_1^{\max}$$

Last statement obtained by using condition 2. In the moment of time $T_2$ control switches. Since $q_2(t) = 0$, $t \geq T_2$, we can consider only $q_1(t)$

$$q_1(t) = q_1(T_2) - \int_{T_2}^t (v(\tau) + \alpha(\tau) - \mu) d\tau,$$

$$q_1(t) \le q_1(T_2) - \varepsilon(t - T_2),$$

$$t \le \frac{q_1(T_2) + \varepsilon T_2}{\varepsilon}.$$

Therefore, we have following estimation for the time $T_1$

$$T_1 \le \frac{q_1(0)}{\varepsilon} + \frac{k}{2}\left(\frac{q_2(0)}{\varepsilon}\right)^2 + T_2 \le \frac{q_1(0)}{\varepsilon} + \frac{q_2(0)}{\varepsilon} + \frac{k}{2}\left(\frac{q_2(0)}{\varepsilon}\right)^2.$$

This completes the proof of Theorem.

**4. Experimental modeling**

Let us illustrate obtained result on the example. Consider the following model:

$$\dot{q}_1(t) = \alpha + q_2(t) - u_1(t),$$

$$\dot{q}_2(t) = v - u_2(t),\ t \ge 0.$$

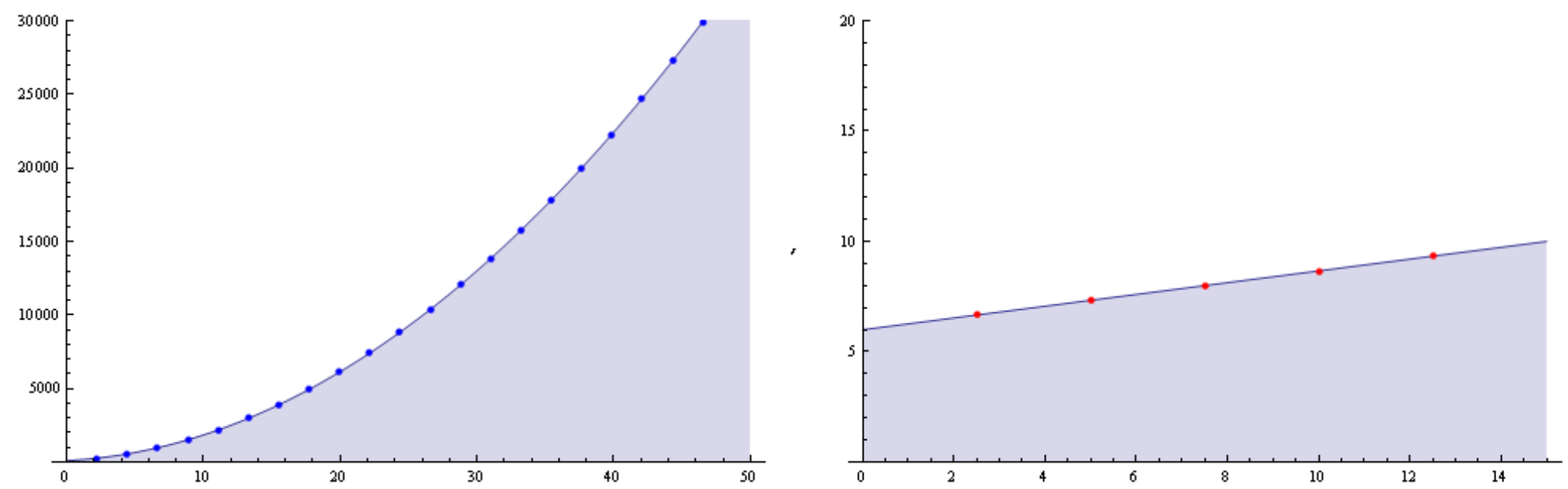

Figure 3. Fluid model without control

On figure 3 we can see dynamic of phase variables $q_1(t)$, $q_2(t)$ without control $u(t)$. Using control defined in previous chapter

$$u(t) = \begin{cases} (\alpha, \mu - \alpha), & t \in [0, T_2] \\ (\mu - v, v), & t \in [T_2, T_1] \end{cases}$$

we can satisfy condition $q(t) = \{0\}$ at the latest at moment of time $T_1$ (Figure 4.).

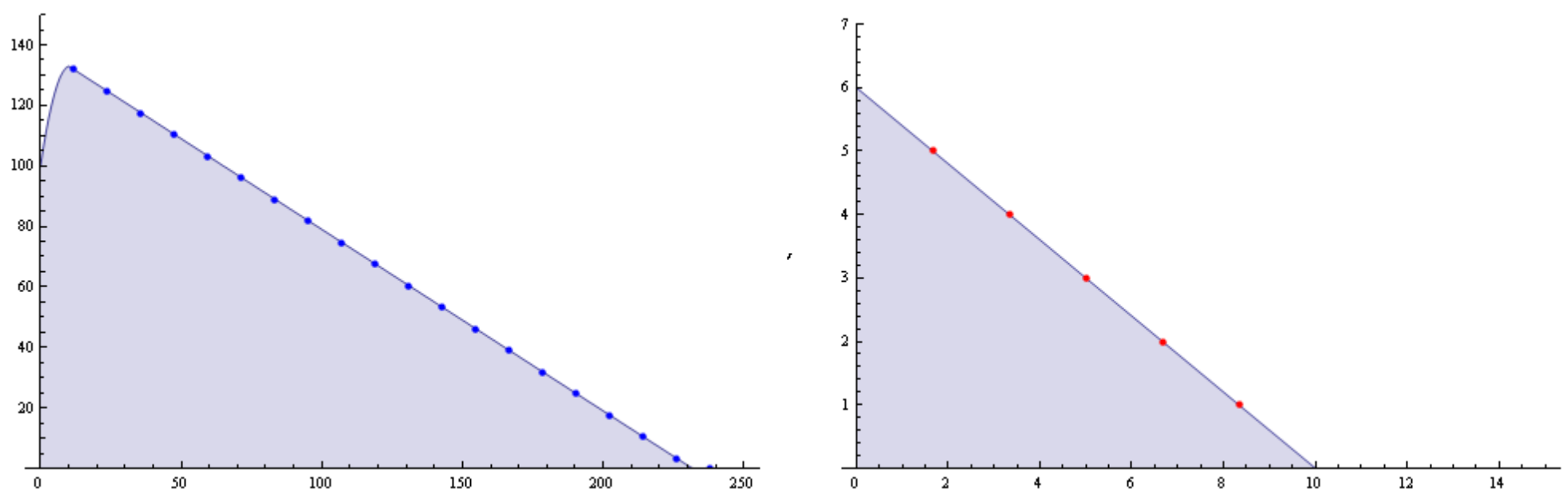

Figure 4. Controlled fluid model

Note that initial state $q_2(0)$ imply growing initial queue state on the value $\frac{k}{2}\frac{(q_2(0))^2}{(\mu - \nu - \alpha_{\max})}$.

## Conclusion

In this work we consider the single server model with counteraction. We extend the classic fluid model using results from the theory of conflict controlled processes. Our new model describes conflict processes on the networks. For that model we have defined two players – an attacker and a defender. Behavior of the system depends on strategies of these players.

We formulated conditions that guarantee possibility of winning of one player. This fact proved by the main theorem. We have illustrated theoretical results on the experimental computational modeling.

## References

1. Chlamtac I., Jain R. Methodology for Building a Simulation Model for Efficient Design and Performance Analysis of Local Area Networks // Simulation, Vol. 42, No. 2, February 1984, P. 57-66.
2. O. Ignatenko, P. Andon. Counteraction to denial of service attacks in Internet: approach concept // Problems of programming, 2-3, 2008, P. 564-574.
3. O. Ignatenko. Denial of service attack in the Internet: agent-based intrusion detection and reaction. http://arxiv.org/abs/0904.4174v1.
4. *Meyn S.* Control Techniques for Complex Networks. – Cambridge University Press, 2007. – 582 p.
5. Bertsekas D. Network optimization: continuous and discrete models. – Athena Scientific, Belmond, 1998. – 270 p.
6. Chikrii A.A. Conflict-controlled Processes, Kluwer Academic Publisher, Boston-London-Dordrecht, 1997, 424 p.
7. M.S. Nikolskii, L.S. Pontryagin's First Direct Method in Differential Games, Izdat. Gos. Univ., Moscow, 1984, 65 p.
8. L.S. Pontryagin, Selected Scientific Papers, Vol. 2, Nauka, Moscow, 1988, 576 p.